\documentclass[11pt]{amsart}
\usepackage{amssymb,amsmath,amsthm,verbatim,mathrsfs}

\newtheorem{thm}{Theorem}
\newtheorem*{thm*}{Theorem}
\newtheorem{prop}[thm]{Proposition}
\newtheorem{defn}{Definition}
\newtheorem{lem}[thm]{Lemma}

\theoremstyle{remark}

\renewcommand{\epsilon}{\varepsilon}

\def\ve{\varepsilon}

\title{The Mabuchi metric and the K\"ahler-Ricci flow.}

\author{Donovan McFeron }
\address{School of Theoretical and Applied Science, Ramapo College of New Jersey, Mahwah, NJ}
\email{dmcferon@ramapo.edu}

\begin{document}

\begin{abstract}
In this paper we show that on a Fano manifold the convergence of the K\"ahler-Ricci flow to a K\"ahler-Einstein metric follows from the integrability of the $L^2$ norm of the Ricci potential for positive time.
\end{abstract}

\maketitle

\section{Introduction}
Let $(M,\omega)$ be a smooth compact Fano manifold of dimension $n$ without boundary, and let the K\"ahler-Ricci flow, which we will always assume  to be volume normalized,  be defined by (see Hamilton~\cite{H1982})
\begin{equation}
\frac{\partial}{\partial t}g_{i\bar j}=-(R_{i\bar j}-g_{i\bar j}), \hspace{5mm} g_{i\bar j}(0)=(g_0)_{i\bar j},
\end{equation}
where $(g_0)_{i\bar j}$ is a given initial K\"ahler metric. We define the Ricci potential $u=u(t)$ by
\begin{equation}
-\partial_i\partial_{\bar j}u=R_{i\bar j}-g_{i\bar j}.
\end{equation}
Under the normalization given by
\begin{equation}
\frac{1}{V}\int_Me^{-u}\omega^n=1, \hspace{5mm} V=\int_M\omega^n,
\end{equation}
it is easy to check that $u$ evolves by
\begin{equation}
\frac{\partial}{\partial t}u=\Delta u+u-b,
\end{equation}
where $b=b(t)$ is the average of $u$ with respect to the measure $e^{-u}\omega^n$.
\begin{equation}
b=\frac{1}{V}\int_Mue^{-u}\omega^n.
\end{equation}
From the work of Perelman~\cite{P2002}(see Sesum-Tian~\cite{ST2008}), along the flow for some constant $C$ depending only on $g_{i\bar j}(0)$,
\begin{equation}\label{C0boundsu}
\|u\|_{C^0}+\|\nabla u\|_{C^0}+\|\Delta u\|_{C^0}\leq C.
\end{equation}
In Theorem~\ref{FiniteMabuchi} (stated below), we assume that the integral over positive time of the $L^2$ norm of the Ricci potential is finite. We will get the sharpest results by choosing the normalization that minimizes the $L^2$ norm, which is the one obtained by taking the difference between the Ricci potential and its average. We denote it by
\begin{equation}
\tilde u=u-a,
\end{equation}
where $a=a(t)$ is the average of $u$ with respect to the measure $\omega^n$.
\begin{equation}
a=\frac{1}{V}\int_Mu\omega^n.
\end{equation}
(Notice that $\int_M\tilde u\omega^n=0$). By direct calculation and integration by parts, $\tilde u$ evolves by
\begin{equation}
\frac{\partial}{\partial t}\tilde u=\Delta \tilde u+\tilde u+\frac{1}{V}\int_M|\nabla\tilde u|^2\omega^n.
\end{equation}
Corresponding to the $C^0$ bounds in Equation~(\ref{C0boundsu}), along the K\"ahler-Ricci flow we have the uniform $C^0$ bounds 
\begin{equation}\label{C0boundstildeu}
\|\tilde u\|_{C^0}+\|\nabla \tilde u\|_{C^0}+\|\Delta \tilde u\|_{C^0}\leq C,
\end{equation}
where $C$ is some positive constant depending only on the initial metric.

The space of K\"ahler potentials can be equipped with a variety of metrics.  For the study of K\"ahler-Einstein metrics, one of the most natural is the Mabuchi metric (see Mabuchi~\cite{M1987}). Another natural metric is the Calabi metric, which is recently studied in connection with the K\"ahler-Ricci flow by Clarke-Rubinstein~\cite{CR2011}. We will only be interested in the length with respect to these metrics, which we define as follows.
\begin{defn}
The length of $u$ with respect to the Mabuchi metric is
\begin{equation}
\int_0^\infty \sqrt{\int_Mu^2_t\omega^n_t}dt=\int_0^\infty\|u\|_{L^2(t)}dt.
\end{equation}
\end{defn}
\begin{defn}
The length of $u$ with respect to the Calabi metric is
\begin{equation}
\int_0^\infty\|\Delta u\|_{L^2(t)}dt.
\end{equation}
\end{defn}

In this paper, we prove the following theorem.

\begin{thm}\label{FiniteMabuchi}
If the length of the normalized Ricci potential, $\tilde u$, with respect to the Mabuchi metric is finite, then the K\"ahler-Ricci flow converges smoothly and exponentially fast to a K\"ahler-Einstein metric.
\end{thm}

The motivation for this theorem came in part from Clarke-Rubinstein~\cite{CR2011}, in which the following theorem is proven.
\begin{thm}[\cite{CR2011} Corollary 6.8]
If the length of the Ricci potential, $u$, with respect to the Calabi metric is finite, then the K\"ahler-Ricci flow converges smoothly and exponentially fast to a K\"ahler-Einstein metric.
\end{thm}

This is a consequence of Theorem~$6.6$ in Clarke-Rubinstein~\cite{CR2011}, which states that if the flow converges with respect to the Calabi metric, then it converges smoothly. It would be interesting to consider if the same result is true for the Mabuchi metric. We will show that due to smoothing properties along the K\"ahler-Ricci flow the length with respect to the Mabuchi metric controls the length with respect to the Calabi metric.  One could proceed with the argument in Clarke-Rubinstein~\cite{CR2011} to show convergence of the K\"ahler Ricci flow. Instead we prove a further smoothing result, which allows us to apply a theorem of Phong-Song-Sturm-Weinkove~\cite{PSSW2009}.

In Section~\ref{sec:Smooth} we make use of two results from the work of Perelman~\cite{P2002}(see Sesum-Tian~\cite{ST2008}). First we use the $C^0$ bounds for $\tilde u$ given in equation~(\ref{C0boundstildeu}) to show the $L^2$ norm of $\tilde u$ controls the $L^2$ norms of $\nabla\tilde u$ and $\Delta\tilde u$.  Second, the monotonicity of Perelman's functional implies the logarithmic Sobolev inequality stated in Proposition~\ref{logSobolev}. Then following the method of Davies~\cite{D1989} (see also Ye~\cite{Y2007}, Zhang~\cite{Z2007} and Cao-Zhang~\cite{CZ2010}), we establish a smoothing lemma, which is an improvement on the following lemma in Phong-Song-Sturm-Weinkove~\cite{PSSW2009}.
\begin{lem}[\cite{PSSW2009} Lemma 1]\label{PSSWsmooth}
There exist positive constants $\delta$ and $K$ depending only on $n$ with the following property. For any $\ve$ with $0<\ve\leq\delta$ and any $t\geq0$, if
\begin{equation*}
\|u\|_{C^0(t)}\leq\ve,
\end{equation*}
then
\begin{equation*}
\|\nabla u\|_{C^0(t+2)}+\|R-n\|_{C^0(t)}\leq K\|u\|_{C^0(t)}.
\end{equation*}
\end{lem}

Finally, in Section~\ref{sec:Main} we prove Theorem~\ref{FiniteMabuchi} by making use of the smoothing results from Section~\ref{sec:Smooth} and by imploying the following lemma from Phong-Song-Sturm-Weinkove~\cite{PSSW2009}.
\begin{lem}[\cite{PSSW2009} Lemma 6]\label{PSSWconvergence}
Assume that the scalar curvature $R(t)$ along the K\"ahler-Ricci flow satisfies
\begin{equation}
\int_0^\infty\|R-n\|_{C^0(t)}dt<\infty.
\end{equation}
Then the metrics $g_{i\bar j}$ converge exponentially fast in $C^\infty$ to a K\"ahler-Einstein metric.
\end{lem}

There have been many other works on the K\"ahler-Ricci flow on Fano manifolds under various assumptions. For the case when a K\"ahler-Einstein metric exists see Tian-Zhu~\cite{TZ2011}. For results which assume various ``stability'' conditions see for example Phong-Sturm~\cite{PS2006}, Sz\'ekelyhidi~\cite{S2010}, and Tosatti~\cite{T2010}.

\subsection*{Acknowledgements}
I would like to thank G\'abor Sz\'ekelyhidi for helpful discussions during the writing of this paper, and Duong H.~Phong for his encouragement and interest in this work. I would also like to thank Yanir Rubinstein for comments on a previous version.

\section{Smoothings}\label{sec:Smooth}

In this section we establish a smoothing lemma on the Ricci potential.  This will be the key tool in the proof of Theorem~\ref{FiniteMabuchi}.
First we state the following elementary lemma.
\begin{lem}\label{ODE}
Suppose that the function $F(t)$ satisfies
\begin{equation}
\frac{\partial}{\partial t}F(t)\leq kF(t),
\end{equation}
for some time independent constant $k$. Then 
\begin{equation}
F(t+1)\leq e^{k}F(t).
\end{equation}
\end{lem}

Next, following the method of Bando~\cite{B1987} (See also Phong-Song-Sturm-Weinkove~\cite{PSSW2009}), we calculate:
\begin{equation}\label{timedervs}
\begin{aligned}
\frac{\partial}{\partial t}\tilde u^2&= \Delta \tilde u^2-2|\nabla \tilde u|^2+2\tilde u^2+\frac{2\tilde u}{V}\int_M|\nabla\tilde u|^2\omega^n,
\\
\frac{\partial}{\partial t}|\nabla \tilde u|^2 &= \Delta|\nabla \tilde u|^2-|\nabla\overline\nabla \tilde u|^2-|\nabla^2\tilde u|^2+|\nabla \tilde u|^2,
\\
\frac{\partial}{\partial t}\Delta \tilde u&= \Delta^2 \tilde u+\Delta \tilde u+|\nabla\overline\nabla \tilde u|^2.
\end{aligned}
\end{equation}
We will rely on these calculations in the proofs of the following lemmas.
\begin{lem}\label{ucontrolsgradu}
There exists a constant $C$ independent of time, such that for all $t\geq 0$,
\begin{equation}
\|\nabla \tilde u\|_{L^2(t+1)}\leq C\|\tilde u\|_{L^2(t)}.
\end{equation}
\end{lem}
\begin{proof}
Let $T\geq0$ be a fixed time. Consider the function $F(t)$ defined by
\begin{equation}
F(t)=\int\left[(t-T)|\nabla \tilde u|^2+D\tilde u^2\right]\omega^n,
\end{equation}
for $t\in[T,T+1]$ and for any time independent constant $D\geq 1/2$. Making use of equations (\ref{C0boundstildeu}) and (\ref{timedervs}) and integration by parts, we calculate
\begin{equation}
\begin{aligned}
\frac{\partial}{\partial t}F(t) 
=&\int_M\left\{\left[(t-T)|\nabla\tilde u|^2+D\tilde u^2\right]\Delta\tilde u\right.
\\
&
+|\nabla\tilde u|^2+(t-T)\left[\Delta|\nabla\tilde u|^2-|\nabla\overline\nabla\tilde u|^2-|\nabla^2\tilde u|^2+|\nabla\tilde u|^2\right]
\\
&\left.
+D\Delta\tilde u^2-2D|\nabla\tilde u|^2+2D\tilde u^2+\frac{2D\tilde u}{V}\int_M|\nabla\tilde u|^2\omega^n
\right\}\omega^n
\\
\leq& \int_M\{\|\Delta\tilde u\|_{C^0}\left[(t-T)|\nabla\tilde u|^2+D\tilde u^2\right]+(t-T)|\nabla\tilde u|^2
\\
&+2D\tilde u^2+(1-2D)|\nabla\tilde u|^2\}\omega^n.
\\
\leq& \int_M(2+\|\Delta\tilde u\|_{C^0})\left[(t-T)|\nabla\tilde u|^2+D\tilde u^2\right]\omega^n\leq kF(t).
\end{aligned}
\end{equation}
where $k$ is a positive constant independent of time. By Lemma~\ref{ODE} we have
\begin{equation}
\int\left(|\nabla\tilde u|_{(T+1)}^2+D\tilde u_{(T+1)}^2\right)\omega_{(T+1)}^n
\leq e^kD\int\tilde u_{(T)}^2\omega_{(T)}^n.
\end{equation}
This completes the proof of Lemma~\ref{ucontrolsgradu}.
\end{proof}

\begin{lem}\label{graducontrolslaplaceu}	
There exists a constant $C$ independent of time, such that for all $t\geq 0$,
\begin{equation}
\|\Delta \tilde u\|_{L^2(t+1)}\leq C\|\nabla \tilde u\|_{L^2(t)}.
\end{equation}
\end{lem}
\begin{proof}
Let $T\geq0$ be a fixed time. Consider the function $F(t)$ defined by
\begin{equation}
F(t)=\int\left[(t-T)|\Delta\tilde u|^2+D|\nabla\tilde u|^2\right]\omega^n,
\end{equation}
for $t\in[T,T+1]$ and for any time independent constant $D\geq n+2\|\Delta\tilde u|_{C^0}$. Making use of equations (\ref{C0boundstildeu}) and (\ref{timedervs}), integration by parts, and the 
estimate $|\Delta\tilde u|^2\leq n|\nabla\overline\nabla\tilde u|^2$, which follows from Cauchy-Schwarz, we calculate
\begin{equation}
\begin{aligned}
\frac{\partial}{\partial t}F(t) 
=&\int_M\left\{\left[(t-T)|\Delta\tilde u|^2+D|\tilde\nabla u|^2\right]\Delta\tilde u\right.
\\
&+|\Delta\tilde u|^2+(t-T)\left[-2|\nabla\Delta\tilde u|^2+2|\Delta\tilde u|^2+2\Delta\tilde u|\nabla\overline\nabla\tilde u|^2\right]
\\
&\left.+D\left[\Delta|\nabla\tilde u|^2-|\nabla\overline\nabla\tilde u|^2-|\nabla^2\tilde u|^2+|\nabla\tilde u|^2\right]\right\}\omega^n
\\
\leq&\int_M\left\{\|\Delta\tilde u\|_{C^0}\left[(t-T)|\Delta\tilde u|^2+D|\tilde\nabla u|^2\right]+2(t-T)|\Delta\tilde u|^2\right.
\\
&\left.+\left(n+2\|\Delta\tilde u\|_{C^0}-D\right)|\nabla\overline\nabla\tilde u|^2+D|\nabla\tilde u|^2\right\}\omega^n
\\
\leq& \int_M(2+\|\Delta\tilde u\|_{C^0})\left[(t-T)|\Delta u|^2+D|\nabla u|^2\right]\omega^n\leq kF(t),
\end{aligned}
\end{equation}
where $k$ is a positive constant independent of time. By Lemma~\ref{ODE} we have
\begin{equation}
\int_M\left(|\Delta\tilde u|_{(T+1)}^2+D|\nabla\tilde u|_{(T+1)}^2\right)\omega_{(T+1)}^n
\leq e^{k}D\int_M |\nabla\tilde u|_{(T)}^2\omega_{(T)}^n.
\end{equation}
This completes the proof of Lemma~\ref{graducontrolslaplaceu}
\end{proof}

The monotonicity of Perelman's functional~\cite{P2002}(see also Ye~\cite{Y2007} and Zhang~\cite{Z2007}), implies the following logarithmic Sobolev inequality.
\begin{prop}\label{logSobolev}
Along the normalized K\"ahler Ricci flow, for all $v\in W^{1,2}(M)$ with $\|v\|_{L^2}=1$ and for any positive constant $A$, there exists a constant $C$ depending only on the initial metric and $A$ such that $\forall\ve\in(0,A)$ 
\begin{equation}\label{eq:logSobolev}
\int_Mv^2\ln v \omega^n \leq \ve\int_M|\nabla v|^2\omega^n-\frac{n}{2}\ln\ve+C.
\end{equation}
\end{prop}

\begin{prop}\label{heatkernel}
Along the normalized K\"ahler Ricci flow, for any non-negative function $f$ satisfying $\frac{\partial}{\partial t}f\leq\Delta f$, there exists a constant $C$ depending only on the initial metric such that for all $t\geq 0$,
\begin{equation}
\|f\|_{C^0(t+1)}\leq C\|f\|_{L^1(t)}.
\end{equation}
\end{prop}
\begin{proof}
Following the method of Davies~\cite{D1989} (see also Ye~\cite{Y2007}, Zhang~\cite{Z2007} and Cao-Zhang~\cite{CZ2010}), for a fixed $T\geq 0$ and $t\in(T,T+1)$ we let 
\begin{equation*}
p(t)=\frac{1}{T+1-t}.
\end{equation*}
Notice that we have the following bounds on $p$.
\begin{equation}\label{pbounds}
1<p<\infty,
\hspace{5mm}
0<\frac{1}{p}<1,
\hspace{5mm}
0<\frac{4(p-1)}{p^2}\leq1.
\end{equation}
By direct computation and integration by parts, we calculate
\begin{equation}
\begin{aligned}
\frac{\partial}{\partial t}\ln\|f\|_{L^p} \leq& 
\frac{p'}{p}\int_M\frac{f^p}{\|f\|_{L^p}^p}\ln f\omega^n
-\frac{p'}{p^2}\ln\int_Mf^p\omega^n
\\
&
-\frac{4(p-1)}{p^2}\int_M\frac{|\nabla (f^{\frac{p}{2}})|^2}{\|f\|_{L^p}^p}\omega^n
-\frac{1}{p}\int_M\frac{Rf^p}{\|f\|_{L^p}^p}\omega^n+\frac{n}{p}
\end{aligned}
\end{equation}
After making the substitution
\begin{equation}
v=\frac{f^{\frac{p}{2}}}{\|f^{\frac{p}{2}}\|_{L^2}}
\end{equation}
and using the bounds in Equation~(\ref{pbounds}), we have
\begin{equation}
\begin{aligned}
\frac{\partial}{\partial t}\ln\|f\|_{L^p} \leq& 
\frac{p'}{p^2}\int_Mv^2\ln v^2\omega^n
-\frac{4(p-1)}{p^2}\int_M|\nabla v|^2\omega^n
-\frac{1}{p}\int_M Rv^2\omega^n+\frac{n}{p}
\\
\leq&
\int_Mv^2\ln v^2\omega^n
-\frac{4(p-1)}{p^2}\int_M|\nabla v|^2\omega^n
+\sup R^-(x,t)+n.
\end{aligned}
\end{equation}
Now we make use of the fact that $\sup R^-(x,t)\leq\sup R^-(x,0)$ and apply the logarithmic Sobolev inequality in Proposition~\ref{logSobolev}
with 
\begin{equation*}
\ve=\frac{4(p-1)}{p^2}=4(T+1-t)(t-T),
\end{equation*}
which yields
\begin{equation}
\begin{aligned}
\frac{\partial}{\partial t}\ln\|f\|_{L^p} \leq& 
-\frac{n}{2}\ln[4(T+1-t)(t-T)]+\sup R^-(x,0)+n+C.
\end{aligned}
\end{equation}
Lastly, integrating from $t=T$ to $t=T+1$, we reach the desired inequality
\begin{equation}
\|f\|_{L^\infty(T+1)}\leq C_1\|f\|_{L^1(T)}
\end{equation}
for all $T\geq 0$ and some positive constant $C_1$ depending only on 
the initial metric.
\end{proof}

\begin{lem}\label{smooth}
Along the normalized K\"ahler Ricci flow, there exists a positive constant $C$ depending only on the initial metric such that for all $t\geq 0$
\begin{equation}
\|\Delta\tilde u\|_{C^0(t+3)}+\|\nabla\tilde u\|_{C^0(t+3)}\leq C\|\tilde u\|_{L^2(t)}.
\end{equation}
\end{lem}
\begin{proof}
Consider the positive function $f$ defined by
\begin{equation}
f(x,t)=e^{-2t}\left[(\Delta\tilde u)^2+D|\nabla\tilde u|^2\right],
\end{equation}
where $D\geq 2\|\Delta\tilde u\|_{C^0}$ is independent of time. Using equation~(\ref{timedervs}) we have
\begin{equation}
\left(\frac{\partial}{\partial t}-\Delta\right)f=e^{-2t}\left[(2\Delta\tilde u-D)|\nabla\overline\nabla\tilde u|^2 -2|\nabla\Delta\tilde u|^2-D|\nabla^2\tilde u|^2\right]\leq0.
\end{equation}
$f$ satisfies the conditions of Proposition~\ref{heatkernel}, so we have
\begin{equation}
\|f\|_{C^0(t+3)}\leq C_1\|f\|_{L^1(t+2)},
\end{equation} 
for some positive constant $C_1$ depending only on the initial metric.
It follows that for some positive constant $C_2$ depending only on the initial metric, 
\begin{equation}
\begin{aligned}
\|\Delta\tilde u\|_{C^0(t+3)}^2+\|\nabla\tilde u\|_{C^0(t+3)}^2\leq
&
C_2\left(\|\Delta\tilde u\|_{L^2(t+2)}^2+\|\nabla\tilde u\|_{L^2(t+2)}^2\right).
\end{aligned}
\end{equation}
We complete the proof by applying Lemma~\ref{ucontrolsgradu} and Lemma~\ref{graducontrolslaplaceu}, which yield
\begin{equation}
\|\Delta\tilde u\|_{C^0(t+3)}^2+\|\nabla\tilde u\|_{C^0(t+3)}^2
\leq C_3\|\nabla\tilde u\|_{L^2(t+1)}^2
\leq C_4\|\tilde u\|_{L^2(t)}^2,
\end{equation}
where $C_3$ and $C_4$ are positive constants depending only on the initial metric.
\end{proof}

\section{Proof of the main theorem}\label{sec:Main}
We finish by proving that Theorem~\ref{FiniteMabuchi} is a consequence of Lemma~\ref{smooth} and Lemma~\ref{PSSWconvergence} (Lemma 6 from Phong-Song-Sturm-Weinkove\cite{PSSW2009}). 
\begin{proof}[Proof of Theorem~\ref{FiniteMabuchi}]
\begin{equation}
\begin{aligned}
\int_0^\infty \|R-n\|_{C^0(t)}dt
=&\int_0^\infty \|\Delta\tilde u\|_{C^0(t)}dt
\\
=&\int_0^3 \|\Delta\tilde u\|_{C^0(t)}dt
+\int_3^\infty \|\Delta\tilde u\|_{C^0(t)}dt
\\
\leq& C_1
+\int_0^\infty \|\Delta\tilde u\|_{C^0(t+3)}dt,
\end{aligned}
\end{equation}
where $C_1$ is independent of time. By Lemma~\ref{smooth} we have
\begin{equation}
\begin{aligned}
\int_0^\infty \|R-n\|_{C^0(t)}dt
\leq& C_1
+C_2\int_0^\infty \|\tilde u\|_{L^2(t)}dt
\leq C.
\end{aligned}
\end{equation}
The final inequality follows from the finiteness of the length of $\tilde u(t)$ with respect to the Mabuchi metric. We are now able to imploy Lemma~\ref{PSSWconvergence} (Lemma $6$ from Phong-Song-Sturm-Weinkove~\cite{PSSW2009}), which concludes that under these conditions the K\"ahler-Ricci flow converges exponentially fast in $C^\infty$ to a K\"ahler-Einstein metric.
\end{proof}

\end{document}